\documentclass[12pt]{article}

\usepackage{amsmath}
\usepackage{amssymb}
\usepackage{amsthm}

\newtheorem{thm}{Theorem}[section]
\newtheorem{prop}[thm]{Proposition}
\newtheorem{defn}[thm]{Definition}
\newtheorem{Exam}[thm]{Example}
\newtheorem{Coro}[thm]{Corollary}
\newtheorem{Lemm}[thm]{Lemma}
\newtheorem{Rema}[thm]{Remark}
\newtheorem{conj}[thm]{Conjecture}
\newcommand{\mult}{\mathrm{mult}}
\newcommand{\ord}{\mathrm{ord}}

\newcommand{\ind}{\mathrm{ind}}
\newcommand{\wgh}{\mathrm{w}}
\newcommand{\wght}{\mathrm{wt}}

\begin{document}
\setcounter{section}{-1}

\begin{center}
{\Large 
Freeness of adjoint linear systems on threefolds 
with non Gorenstein $\mathbb{Q}$-factorial 
terminal singularities or 
some quotient singularities }\\
\end{center}
\begin{center}
 Nobuyuki Kakimi
\end{center}
\begin{flushleft}
 Department of Mathematical Sciences, University of Tokyo, Komaba, Meguro, 
 Tokyo 153, Japan (e-mail:kakimi@318uo.ms.u-tokyo.ac.jp)  
\end{flushleft}
\begin{flushleft}
Abstract. 
We define and calculate the weighted multiplicities of 
non Gorenstein terminal singularities on
threefolds and  some quotient singularities.  
As an application, we improve freeness conditions on threefolds.  
\end{flushleft}
\section{Introduction}
\hspace{1.5 em} 
Let $X$ be a normal projective variety of dimension $3$, 
$x_0 \in X$ a point, and $L$ an ample $\mathbb{Q}$-Cartier divisor 
such that $K_X + L$ is a Cartier divisor at $x_0$. 
Assume that $\sqrt[p]{ L^p \cdot W} \geq {\sigma}_p$ 
for any subvariety $W$ of dimension $p$ which contains $x_0$ 
for some fixed real numbers. 
Then we shall prove the following theorems: 
\begin{thm}[{=Theorem 4.1}] 
Assume $x_0 \in X$ is a quotient singular point 
of type $({1}/{r},{a}/{r},{b}/{r})$
such that an integer $r >0 $, $(r, a) = 1$, and $(r,b) = 1$.
Let ${\sigma}_1 \geq {3}/{r}$,
${\sigma}_2 \geq {3}/{\sqrt{r}}$, and 
${\sigma}_3 >  {3}/{\sqrt[3]{r}}$. 
Then $| K_X + L |$ is free at $x_0$. 
\end{thm} 
\begin{thm}[{=Theorem 4.4}]
Assume $x_0 \in X$ is a nonhypersurface and not quotient 
$\mathbb{Q}$-factorial terminal singular point of $\ind_{x_0} X =r > 1$. 
Let ${\sigma}_1 \geq 2/r$, ${\sigma}_2 \geq 2\sqrt{2/r}$, 
and ${\sigma}_3 > 2\sqrt[3]{2/r}$. 
Then $| K_X + L |$ is free at $x_0$. 
\end{thm}
\begin{Coro}[{= Corollary 4.5}] 
Assume $X$ have only $\mathbb{Q}$-factorial terminal singularities and 
$x_0 \in X$ a point of $\ind_{x_0} X=r > 0$. 
Let ${\sigma}_1 \geq {3}/{r}$,
${\sigma}_2 \geq {3}/{\sqrt{r}}$, and 
${\sigma}_3 >  {3}/{\sqrt[3]{r}}$. 
Then $| K_X + L |$ is free at $x_0$. 
\end{Coro}
The above results are generalization of the following conjecture 
by T. Fujita: 
\begin{conj}
For a smooth projective variety $X$ and 
an ample divisor $L$ on $X$, the linear system $|K_X + mL|$ is free 
if $m \geq dim X + 1$. 
\end{conj}
A strong version of Fujita's freeness conjecture is the following: 
\begin{conj}
Let $X$ be a normal projective variety of dimension $n$, 
$x_0 \in X$ a smooth point, and $L$ an ample Cartier divisor.
Assume that there exist positive numbers ${\sigma}_p$ for $p=1,2, \cdots, n$
which satisfy the following conditions: \\ 
$(1)$ $\sqrt[p]{L^p \cdot W} \geq {\sigma}_p$ for 
any subvariety $W$ of dimension $p$ which contains $x_0$, \\
$(2)$ ${\sigma}_p \geq n$ for all $p$ and ${\sigma}_n > n$. \\
Then $|K_X+L|$ is free at $x_0$. 
\end{conj}
We list the cases where the conjectures are already proved. 
For smooth complex algebraic surface,  
Reider [Rdr] proved the strong version of Fujita's freeness conjecture 
by applying Bogomolov's instability 
theorem to study adjoint series on surfaces. 
For a projective normal surface, Ein and Lazarsfeld [EL], 
Matsushita [Mat], Kawachi [KM][Kwc], and Ma\c{s}ek [Ma] extended 
the result of Reider [Rdr] to singular cases. 
Langer [La1][La2] obtained 
the best estimates for a normal surface 
by applying a rank $2$ reflexive sheaf. 
For Fujita's freeness conjecture, 
it is quite hard in dimension three proved by Ein and Lazarsfeld [EL]. 
The lectures of Lazarsfeld [L] provided a very good introduction. 
Kawamata[Ka5] proved in dimension four case. 
For the strong version of Fujita's freeness conjecture, 
Fujita [F] proved that, if ${\sigma}_1 \geq 3$, ${\sigma}_2 \geq \sqrt{7}$, 
and $\sigma_3 > \sqrt[3]{51}$, then $|K_X + L|$ is free at $x_0$. 
Kawamata [Ka5] proved the following:  
\begin{thm}[{[Ka5, Theorem 3.1]}]
Let $X$ be a normal projective 
variety of dimension $3$, 
$x_0 \in X$ a smooth point, and $L$ an ample Cartier divisor.
Assume that there exist positive numbers ${\sigma}_p$ for $p=1, 2, 3$
which satisfy the following conditions: \\ 
$(1)$ $\sqrt[p]{L^p \cdot W} \geq {\sigma}_p$ for 
any subvariety $W$ of dimension $p$ which contains $x_0$, \\ 
$(2)$ ${\sigma}_1 \geq 3$, ${\sigma}_2 \geq 3$, and ${\sigma}_3> 3$. \\ 
Then $|K_X+L|$ is free at $x_0$. 
\end{thm}
For a projective variety $X$ of dimension $3$ with some singularities, 
Oguiso and Peternell [OP] proved that, 
with only $\mathbb{Q}$-factorial Gorenstein terminal 
$(\textrm{resp.}\ \textrm{\ canonical\ })$ singularities 
and an ample divisor $L$ on $X$, 
the linear system $|K_X + mL|$ is free 
if $m \geq 5 \ (\textrm{resp.}\ m \geq 7)$. 
Ein, Lazarsfeld and Ma\c{s}ek [ELM], 
and Matsushita [Mat] extended some of the results 
of Ein and Lazarsfeld [EL] to projective threefolds with 
terminal singularities. 
We [K1] extended the result of Kawamata [Ka5] 
to normal projective threefolds 
with terminal Gorenstein singularities or some quotient singularities. 
Our freeness conditions in terminal Gorenstein singularities 
or quotient singular points of type $1/r(1,1,1)$ 
are better than in smooth case. 
We noticed that our proof [K1, Theorem 3.8] of canonical 
and not terminal singular point case is wrong. 
Note that Lee [L1] [L2] also obtained some results on 
only Gorenstein canonical singularities independently. \\

The idea of this paper is to define the weighted multiplicities 
of normal singularities and control the constants $\sigma_p$ 
using this invariants. 
We calculate the weighted multiplicities of some quotient singularities 
and non Gorenstein terminal singularities on threefolds. 
\begin{thm}[{= Theorem 2.5}]
Let $(X, x_0) \cong \mathbb{C}^n/{\mathbb{Z}_r(1,a_1, \cdots a_{n-1})}$
such that an integer $r > 0$, $(r, a_1) = 1$, 
and integers $a_j$ $( 0 \leq a_j < r )$ for $1 \leq j \leq n-1$. 
Let $l := \min \{ i\ |\ a_{j}i \equiv i 
(\textrm{\ mod\ } r )\ ( 1 \leq j \leq n - 1) 
\textrm{\ for\ } 0 <  i \leq r \}$. 
Let ${\mu}: Y \rightarrow X$ be 
the weighted blow up of $X$ at $x_0$ such that 
$wt(x_0, x_1, \cdots, x_{n-1})= (l/r, l/r, \cdots, l/r)$ 
with the exceptional divisor $E$ of $\mu$. 
Then we have 
\[ {\wgh}\textrm{-}{\mult}_{{\mu}:x_0} X = r^{n-1}/{l}^n.\]  
\end{thm} 
\begin{thm}[{=Theorem 2.7}]
Let $(X, x_0)$ be a $3$ folds nonhypersurface and not quotient 
terminal singular point of 
$\ind_{x_0} X = r >1$ over $\mathbb{C}$ and $\mu:Y \rightarrow X$ 
the weighted blow up with the weights 
${\wght}(x, y, z, u) = (1, 1, 1, 1)$ 
with the exceptional divisor ${E}$ of ${\mu}$
such that $K_Y = {\mu}^{*} K_X + E$. 
Then 
\[{\wgh}\textrm{-}{\mult}_{\mu:x_0} X= 2/r. \]   
\end{thm}
By applying the weighted multiplicities, 
we improve freeness conditions on threefolds (Theorems 4.1, 4.4). 
Our freeness conditions are better than in smooth case. 
We give explicit examples (Examples 4.2, 4.3) 
which show the estimate in Theorem 4.1 is 
best.
Moreover, we obtain an alternative proof (Theorem 5.1) of 
a theorem of Langer for log terminal singular points of only $A$-type [La2]. 
The method of [La2] is to use a reflexive sheaf of rank $2$, 
which is completely different for ours. 
We also give explicit examples (Example 5.2, 5.3)
for Theorem 5.1 which shows the estimate is best. \\
Our proof of freeness is very similar to the one for the smooth case 
given in [Ka5]. However this involves more careful and detailed analysis 
of weighted blow up and the weighted multiplicities. \\
\\
{\bf Acknowledgment}:
The author would like to express his thanks to 
Professor Yujiro Kawamata for his advice and warm encouragement. 
He has benefited from discussions with 
Doctor Masayuki Kawakita, Doctor Daisuke Matsushita, and 
Doctor Hiroshi Sato. 
He also would like to express his thanks to 
Doctor Adrian Langer for his pointing out him his ICTP preprint 
and his kindly comments. 
\section{Preliminaries}
\hspace{1.5 em}
Most of the results of freeness in this paper are the applications of 
the following vanishing theorem:
\begin{thm}[{ [Ka1] [V] }] 
Let $X$ be a smooth projective variety and 
$D$ a ${\mathbb Q}$-divisor.
Assume that $D$ is nef and big, and that the support
of the difference $ {}^{\lceil} {D} { }^{\rceil} - D $ is a normal crossing
divisor. Then ${H}^p (X, K_X + {}^{\lceil} {D} { }^{\rceil}) = 0 $ 
for $p > 0$.
\end{thm}
We recall notation of [Ka5] \ (cf [KMM]).
 \begin{defn} \normalfont
 Let $X$ be a normal variety and $D =\Sigma_i {d}_i {D}_i$ 
 an effective ${\mathbb Q}$-divisor 
 such that $K_X + D$ is ${\mathbb Q}$-Cartier.
 If $\mu:Y \rightarrow X$ is an embedded resolution of the pair $(X,D)$,
 then we can write 
 \[ K_Y + F = {\mu}^* (K_X + D) \]
 with $F = {\mu}_{*}^{-1}D + {\Sigma}_j {e}_j {E}_j$ 
 for the exceptional divisors $E_j$.

 The pair $( X,D )$ is said to have only 
 \textit{log canonical singularities} $(LC) \\
 (\mathrm{resp. } \mathit{kawamata\ log\ terminal\ singularities} (KLT))$ 
 if ${d}_i \leq 1 (\mathrm{resp.} < 1)$
 for all $i$ and  ${e}_j \leq 1(\mathrm{resp.} <1)$ for all $j$.

 A subvariety $W$ of $X$ is said to be a \textit{center of
 log canonical singularities} for the pair $( X,D )$,
 if there is a birational morphism from a normal variety 
 $\mu:Y \rightarrow X$ and a prime divisor $E$ on $Y$ 
 with the coefficient $e \geq 1$
 such that ${\mu}(E) = W$. 
 The set of all the centers of log canonical singularities is denoted
 by $CLC( X, D )$.
 For a point $ x_0 \in X$, we define
 $CLC ( X,x_0,D ) = \{ W \in CLC( X, D ) ; x_0 \in W \}$. 
 \end{defn} 
We shall use the following 
proved by Kawamata [Ka5] [Ka6].
Then we can control the singularities of 
the minimal center of log canonical singularities
and can replace the minimal center of log canonical
singularities by a smaller subvariety.
\begin{prop} [{[Ka5 1.5,1.6,1.9] [Ka6]}]
Let $X$ be a normal variety and  $D$ an effective 
${\mathbb Q}$-Cartier divisor such that $K_X + D $ is ${\mathbb Q}$-Cartier.
Assume that  $X$ is $KLT$ and $(X,D)$ is $LC$.
If $W_1,W_2 \in CLC(X,D)$ and $W$ is 
an irreducible component of $W_1 \cap W_2$,
then $W \in CLC(X,D) $.
If $(X,D)$ is not  $KLT$ at a point $x_0 \in X $,
then there exists the unique minimal element $W_0$ of $CLC(X,x_0,D)$.
Moreover, $W_0$ is normal at $x_0$.   
Then $W_0$ has at most a rational singularity at $x_0$.
\end{prop}
\begin{prop}[{[Ka5 1.10]}] 
Let $x_0 \in X$, $D$ and $W_0$ be as in Proposition $1.3$.
Let $D_1$ and $D_2$ be effective ${\mathbb Q}$-Cartier divisors on $X$ 
whose supports do not contain $W_0$ and which  induce the same
${\mathbb Q}$-Cartier divisor on $W_0$. 
Assume that $(X,D+D_1)$ is $LC$ at $x_0$ and
there exists an element of $CLC(X, x_0, D + {D}_1)$ which is properly
contained in $W_0$. 
Then the similar statement holds for the pair $(X, D + {D}_2)$.
\end{prop} 
We shall need the following Mori's classification theorem 
of terminal singularities in dimension $3$. 
\begin{thm}[{[M]}]
Let $0 \in X$ be a $3$-fold terminal nonhypersurface singular point  
over $\mathbb{C}$.
Then $0 \in X$ is isomorphic to a singularity 
described by the following list:\\
$(1)cA/r, 
\{x y + f(z, {u}^r)=0, f \in \mathbb{C}\{ z, {u}^r \}, (r,a)=1 \}
\subset \mathbb{C}^4/\mathbb{Z}_r(a, r-a, r, 1)$, \\
$(2)cAx/4, \{ x^2+y^2+f(z, u^2)=0, f \in \mathbb{C}\{ z, u^2 \} \}
\subset \mathbb{C}^4/\mathbb{Z}_4(1, 3, 2, 1)$, \\
$(3)cAx/2, \{ x^2 + y^2 +f(z, u)=0,  
f \in (z, u)^4\mathbb{C}\{ z, u \} \}
\subset \mathbb{C}^4/\mathbb{Z}_2(1, 2, 1, 1)$, \\
$(4)cD/2, \{u^2+z^3+x y z + f(x, y)=0, f \in (x, y)^4\}, 
\textrm{or } \{u^2+x y z + {z}^n + f(x, y)=0, f \in (x, y)^4, n \geq 4\}, 
\{u^2+y^2 z + {z}^n + f(x,y)=0, f\in(x, y)^4, n \geq 3 \}
\subset \mathbb{C}^4 / \mathbb{Z}_2(1, 1, 2, 1)$, \\
$(5)cD/3, u^2+x^3+y^3+z^3=0, \textrm{or } 
\{u^2+x^3+y z^2+f(x, y, z)=0, f\in (x, y, z)^4\}, 
\{u^2+x^3+y^3+f(x, y, z)=0, f\in (x, y, z)^4 \}
\subset \mathbb{C}^4 / \mathbb{Z}_3(1, 2, 2, 3)$, \\
$(6)cE/2, \{ u^2+x^3 +g(y, z) x+ h(y, z)=0, g, h 
\in \mathbb{C}\{y, z\}, g, h \in (y, z)^4 \}
\subset \mathbb{C}^4 /\mathbb{Z}_2(2, 1, 1, 1)$. \\
The equations have to satisfy $2$ obvious conditions:
$1$. The equations define a terminal hypersurface singularity.
$2$. The equations are ${\mathbb Z}_n$-equivariant. 
$(\textrm{In fact\ }{\mathbb{Z}_n}\textrm{-invariant,\ except\ for\ }cAx/4.)$
\end{thm}
\section{Definition and Calculation of the weighted multiplicities}
\hspace{1.5 em}
We first define the new following notions which 
we derive from the multiplicity of a point on a normal variety $X$ and 
the multiplicity of an effective $\mathbb{Q}$-Cartier divisor $D$ 
on $X$ at a point:
\begin{defn} \normalfont
Let $X$ be a normal variety of dimension $n$, $x_0$ 
a point of $X$, $\mu: Y \rightarrow X$ a weighted blow up 
at $x_0$ with exceptional divisors $E$, 
$W \subset X$ the subvariety of dimension $p$ 
such that $W$ is normal at $x_0$, $\bar{W}$ the strict transform of $W$, 
and $\bar{D}_{\bar{W}}$ on $\bar{W}$ 
the strict transform of 
an effective $\mathbb{Q}$-Cartier divisor $D_W$ on $W$.  
The \textit{weighted multiplicity} of $W$ at $x_0$ for $\mu$ 
$({\wgh}\textrm{-}\mult_{\mu:x_0} W)$ 
is that
\[ \dim \frac{O_{W, x_0}}{{\mu}_{*}O_{\bar{W}}(-hE|_{\bar{W}})} 
={\wgh}\textrm{-}\mult_{\mu:x_0} W \cdot \frac{h^p}{p!}
+ \textrm{lower\ term\ in\ } h. \]
The \textit{weighted order} of $D_W$ on $W$ at $x_0$ for $\mu$ 
$({\wgh}\textrm{-}{\ord}_{\mu:x_0} D_W)$ 
is that 
\[ {\mu}^{*} (D_W) =
 \bar{ D_{W} } + {\wgh}\textrm{-}{\ord}_{\mu:x_0} D_W \cdot E|_{\bar{W}}. \]
\end{defn}
The weighted multiplicity at a nonsingular point 
for ordinary blow up at the point is 
the same as usual multiplicity at the point. 
Also the weighted order of an effective $\mathbb{Q}$-Cartier 
divisor at a nonsingular point for ordinary blow up at the point is 
the same as usual multiplicity of the divisor at the point. 
We have the following of 
the weighted multiplicity and the weighted order 
at a nonsingular point: 
\begin{Exam} \normalfont 
Let $X = {\mathbb{C}}^3$ and 
$\mu:Y \rightarrow X$ the weighted blow up at $x_0$ 
of the weight $(x,y,z)=(1,a,b)$ for $a,b \in \mathbb{Z}_{> 0}$ 
with the exceptional divisor $E$ 
and 
$D$ an effective $\mathbb{Q}$-Cartier divisor on $X$. 
We consider that 
\[
\frac{ {\mu}_{*}O_{Y}(-hE) }{{\mu}_{*}O_{Y}(-(h+1)E)} 
\ni \begin{array}{@{\,}ll} 
       x^s y^t (y z)^l \textrm{\ for\ } 
       s + t a + l(a+b) =h \textrm{\ or\ } \\
       x^s (y z)^l z^u \textrm{\ for\ } 
       s+ l(a+b)+ u b = h.
    \end{array} \] 
Then we have 
\[ \dim \frac{{\mu}_{*}O_{Y}(-hE)}{{\mu}_{*}O_{Y}(-(h+1)E)} 
= \sum_{l=0}^{ {}_{\lfloor}{h/a+b}_{\rfloor} } 
\{ {\frac{}{} }_{\lfloor} { \frac{h-l(a+b)}{a} }_{\rfloor} + 
{\frac{}{} }_{\lfloor} {\frac{h-l(a+b)}{b} }_{\rfloor} + 1 \}. \]
Then 
\[{\wgh}\textrm{-}{\mult}_{\mu:x_0} X= {1}/{ab}.\]
Moreover, let $F = \sum_{i, j, k} c_{i j k} {x}^i {y}^j {z}^k$ 
be an equation of $D$. 
Then \[ {\wgh}\textrm{-}{\ord}_{\mu:x_0} D 
= \min \{ i + a j + b k ;\  c_{i j k} \neq 0 \}. \]
\end{Exam}
\begin{Exam} \normalfont 
Let $(X,x_0)$ be a $3$ fold terminal singular point of 
${\ind}_{x_0} X =r >1$ over $\mathbb{C}$ 
and $\mu:Y \rightarrow X$ the weighted blow up 
by Kawamata [Ka3] such that 
$K_Y = {\mu}^{*} K_X + (1/r) E$. 
Then,  
\[ {\wgh}\textrm{-}{\mult}_{\mu:x_0} X \leq r^2/(r-1). \]  
Moreover, 
let $(X, x_0) \cong (\mathbb{C}^3/\mathbb{Z}_r(a, r-a, 1), 0)$ or 
$(x y + f(z, {u}^r)=0 \subset \mathbb{C}^4 /\mathbb{Z}_r (a, r-a, r, 1), 0)$ 
for $(r, a) = 1$. Then, 
\[  {\wgh}\textrm{-}{\mult}_{\mu:x_0}X \leq r^2/ a (r-a). \] 
\end{Exam}
For the weighted order, Kawamata [Ka4] or Koll{\'a}r [Ko 8. 13] defined 
the similar notion: 
\begin{Rema} \normalfont
Let $M$ be an effective Weil divisor on 
$( X, x_0 ) = \mathbb{C}^3 / {\mathbb{Z}}_r ( a, - a, 1 )$ 
with $0 < a < r$ and $( r, a ) = 1$, and $\mu : Y \rightarrow X$ 
the weighted blow up at $x_0$ with the weight 
${\wght}(x,y,z)=(a/r, 1 - a/r, 1/r)$ with the exceptional divisor $E$. 
The \textit{weighted multiplicity} of $M$ at $x_0$ 
$({\wgh}-{\mult}_{x_0} M )$ defined by Kawamata [Ka4] is that 
\[ {\wgh}-{\mult}_{x_0} M = 
r \cdot {\wgh}\textrm{-}{\ord}_{\mu:x_0} M. \]
Let $f$ be a holomorphic function near $0 \in \mathbb{C}^n$. 
Assign rational weights $w(x_i)$ to the variables. 
The \textit{weighted multiplicity} of $f$ defined by Koll{\'a}r [Ko 8.13] 
$(w(f))$ is the lowest weight of monomials occurring in $f$.
\end{Rema}  
We calculate the weighted multiplicities in four cases. 
First we calculate the weighted multiplicities  
of some quotient singularities for the weighted blow up at the points. 
\begin{thm}
Let $(X,x_0) \cong \mathbb{C}^n/{\mathbb{Z}_r(1,a_1, \cdots a_{n-1})}$
such that an integer $r > 0$, $(a_1, r)=1$, 
and integers $a_j$ $(0 \leq a_j < r)$ for $1 \leq j \leq n-1$. 
Let $l := \min \{ i\ |\ a_{j}i \equiv i ( \textrm{\ mod\ } r)\ 
( 1 \leq j \leq n-1 ) 
\textrm{\ for\ } 0< i \leq r \}$. 
Let ${\mu}: Y \rightarrow X$ be 
the weighted blow up of $X$ at $x_0$ such that 
$wt(x_0, x_1, \cdots, x_{n-1})= (l/r, l/r, \cdots, l/r)$
with the exceptional divisor $E$ of $\mu$. 
Then we have 
\[ {\wgh}\textrm{-}{\mult}_{{\mu}:x_0} X = {r^{n-1}}/{{l}^n}.\]  
\end{thm}
\begin{proof} 
we consider 
\[ \frac{ O_{X, x_0} }{{\mu}_{*}O_{Y}(-h E)} 
\ni {x_0}^{b_0} {x_1}^{b_1} \cdots {x_{n-1}}^{b_{n-1}} \]
for $h > l( b_0 + b_1 + \cdots + b_{n-1})/r  \in \mathbb{Z}$ and
$b_0 + b_1 a_1 + b_2 a_1 + \cdots + b_{n-1} a_{n-1} \in r \mathbb{Z}$.  

We fix $k_j$ for $1 \leq j \leq n-1$. 
We consider ${b_j}_{k_j}= {m_j}_{k_j} r + k_j$ for 
$0 \leq k_j \leq r-1 $ and ${m_j}_{k_j} \geq 0$.
Then we have 
$b_0 = m r - \overline{ a_1 k_1 + a_2 k_2 + \cdots + a_{n-1} k_{n-1} }$. \\
Then we have 
$l( k_1 + k_2 + \cdots + k_{n-1} 
- \overline{ a_1 k_1 + a_2 k_2 + \cdots + a_{n-1} k_{n-1}}) / r 
\in \mathbb{Z}$ and 
$m + {m_1}_{k_1} + {m_2}_{k_2} + \cdots + {m_{n-1}}_{k_{n-1}} + 
(k_1 + k_2 + \cdots k_{n-1} \\
- \overline{ a_1 k_1 + a_2 k_2 + \cdots + a_{n-1} k_{n-1} })/r \leq (h-1)/l$. 
Then we have 
\[ - r < (k_1 + k_2 + \cdots + k_{n-1} 
- \overline{ a_1 k_1 + a_2 k_2 + \cdots + a_{n-1} k_{n-1} } \leq r(n-1). \]
Then 
\[ r^{n-1} \frac{({}_{\lfloor} {(h-1)/l}_{\rfloor} + 2)!}
{n! ({}_{\lfloor} {(h-1)/l}_{\rfloor} + 2 - n)!} 
\leq \frac{ O_{X, x_0} }{{{\mu}}_{*}O_{Y}(-h E)} 
\leq r^{n-1} \frac{({}_{\lfloor} {(h-1)/l}_{\rfloor} + 2 + n)!}
{n! ({}_{\lfloor} {(h-1)/l}_{\rfloor} + 2)!}. \]  
Hence 
\[{\wgh}\textrm{-}{\mult}_{{\mu}:x_0} X = {r^{n-1}}/{{l}^n}. \]
\end{proof}
Second we calculate the weighted multiplicities of subvariety of 
$\mathbb{C}^3 / {\mathbb{Z}_r(1, a, b)}$. 
\begin{Lemm}
Let $(X,x_0) \cong \mathbb{C}^3/{\mathbb{Z}_r(1, a, b)}$
such that an integer $r > 0$, $(a, r)=1$, and $(b, r)=1$. 
Let $l := \min \{ i\ |\ i \equiv ai \equiv bi ( \textrm{\ mod\ } r)\ 
\textrm{\ for\ } 0< i \leq r \}$. 
Let ${\mu}: Y \rightarrow X$ be 
the weighted blow up of $X$ at $x_0$ such that 
$wt(x, y, z)= (l/r, l/r,l/r)$
with the exceptional divisor $E$ of $\mu$. 
Let $C=\mathbb{C}/{\mathbb{Z}_r(1)}$, 
$S_1 =\mathbb{C}^2/{\mathbb{Z}_r(1, a)}$, 
$S_2= \{ xy + z^n = 0 \subset \mathbb{C}^3/{\mathbb{Z}_r(1, a, b)}\}$, 
and $S_3= \{ x^2 + f(y, z) = 0 \subset \mathbb{C}^3/{\mathbb{Z}_r(1, a, b)}\}$.
Then we have 
\[ {\wgh}\textrm{-}{\mult}_{{\mu}:x_0} C = 1 / l, 
{\wgh}\textrm{-}{\mult}_{{\mu}:x_0} S_1 = r/{l}^2, \]  
\[\textrm{\ and\ } {\wgh}\textrm{-}{\mult}_{{\mu}:x_0} S_i 
= 2r / {l}^2 \textrm{\ for\ } i=2 \textrm{\ or\ } 3, \]
\end{Lemm}
\begin{proof}
Let 
$\bar{C}$ be the strict transform of $C$, 
$\bar{S_i}$ the strict transform of $S_i$ for $1 \leq i \leq 3$, 
and $e$ such that $ae \equiv 1 (\textrm{\ mod\ } r)$. \\
In $C$ case, we consider that, 
\[
\frac{ O_{C, x_0} }{ {\mu}_{*}O_{\bar{C}}(-h E |_{\bar{C}})} 
\ni x^s\textrm{\ for\ } 
h  > s l /r \in \mathbb{Z} \textrm{\ and\ } 
s \in r \mathbb{Z}. \]    
As in Theorem 2.5, we have 
${\wgh}\textrm{-}{\mult}_{{\mu}:x_0} C= 1 /l$. \\ 
In $S_1$ case, we consider that, 
\[
\frac{ O_{S_1, x_0} }{ {\mu}_{*}O_{\bar{S_1}}(-h E |_{\bar{S_1}})} 
\ni x^s y^u \textrm{\ for\ } 
h  > (s + u) l /r \in \mathbb{Z} \textrm{\ and\ } 
s + a u \in r \mathbb{Z}. \]    
As in Theorem 2.5, 
${\wgh}\textrm{-}{\mult}_{{\mu}:x_0} S_1 = r/{l}^2$. \\
In $S_2$ case, we consider that, 
\[
\frac{ O_{S_2, x_0} }{ {\mu}_{*}O_{\bar{S_2}}(-h E |_{\bar{S_2}})} 
\ni \begin{array}{@{\,}ll} 
       x^s z^u \textrm{\ for\ } 
       h  > (s + u) l /r \in \mathbb{Z} \textrm{\ and\ } 
       s + b u \in r \mathbb{Z}
        \textrm{\ or\ } \\
       y^t z^u \textrm{\ for\ } 
       h > (t + u) l /r \in \mathbb{Z} \textrm{\ and\ } 
       t + \overline{e b} u \in r \mathbb{Z}. 
    \end{array} \]  
As in Theorem 2.5, 
${\wgh}\textrm{-}{\mult}_{{\mu}:x_0} S_2 = 2r/{l}^2$. \\
In $S_3$ case, we consider that 
\[ \frac{ O_{S_3, x_0} }{{\mu}_{*}O_{\bar{S_3}}(-h E |_{\bar{S_3}})} 
\ni \begin{array}{@{\,}ll} 
       y^t z^u \textrm{\ for\ } 
       h >(t + u) l/ r \in  \mathbb{Z} \textrm{\ and\ }
       t + \overline{e b} u \in r \mathbb{Z}
        \textrm{\ or\ } \\
        x y^{t} z^{u} \textrm{\ for\ } 
       h >(1 + t + u) l / r \in \mathbb{Z} \textrm{\ and\ }
       e + t +\overline{e b} u \in r \mathbb{Z}.
    \end{array} \]   
We fix $j$ for $0 \leq j < r$.    
We consider $u = m r + j$. \\
We have 
$t = n_1 r - \overline{e b j}$ or $n_2 r - \overline{ e b j } -e$.
Then $(m+ n_1) r + j - \overline{e b j}$ or 
$(m + n_2) r + j - \overline{ e b j }- e +1 \leq r (h-1) / l$. 
Then we have $-r < j - \overline{e b j} < r$ and 
$-2r < j - \overline{e b j}- e + 1 < r$. 
Then
\[ 2r \frac{({}_{\lfloor} {( h - 1) / l}_{\rfloor})!}
{({}_{\lfloor} {( h - 1) / l}_{\rfloor} -2)! 2!}
 < \frac{ O_{S, x_0} }{{\mu}_{*}O_{\bar{S}}(-h E|_{\bar{S}})} < 
2r \frac{({}_{\lfloor} {( h - 1) / l}_{\rfloor} +3)!}
{({}_{\lfloor} {( h - 1) / l}_{\rfloor} +1)!2!}. \]
Hence we have 
${\wgh}\textrm{-}\mult_{\mu:x_0} S= 2r/{l}^2$. 
\end{proof}
Third we calculate the weighted multiplicities of $3$ folds 
nonhypersurface and not quotient 
terminal singular point of $\ind_{x_0} X = r > 1$ 
for a weighted blow up. 
\begin{thm}
Let $(X, x_0)$ be a $3$ folds nonhypersurface and 
not quotient terminal singular point of 
$\ind_{x_0} X = r >0$ over $\mathbb{C}$ and $\mu:Y \rightarrow X$ 
the weighted blow up with the weights 
${\wght}(x, y, z, u) = (1, 1, 1, 1)$  
with the exceptional divisor ${E}$ of ${\mu}$
such that $K_Y = {\mu}^{*} K_X + E$.  
Then 
\[{\wgh}\textrm{-}{\mult}_{\mu:x_0} X= {2}/{r}. \]   
\end{thm}
\begin{proof}
Case$(1)$ 
Let $(X,x_0) \cong ( xy +f(z,u^r)=0 
\subset \mathbb{C}^4/\mathbb{Z}_r(a,r-a,r,1), 0)$ for $(r,a)=1$.
We consider that, for $x y =- f(z, u^r)\ (z \notin f(z, w^r))$, 
\[ \frac{{\mu}_{*}O_{Y}(-h {E})}{{\mu}_{*}O_{Y}(-(h + 1){E})} 
\ni \begin{array}{@{\,}ll} 
       {x}^s {z}^l {u}^m \textrm{\ for\ } 
       s + l + m = h \textrm{\ and\ }
       a s + l r + m \in r \mathbb{Z}
       \textrm{\ or\ } \\
       {y}^t {z}^l {u}^m \textrm{\ for\ } 
       t + l + m = h \textrm{\ and\ }
       (r - a) t +l r + m \in r \mathbb{Z}.
    \end{array} \]
Then, as Theorem 2.5, 
${\wgh}\textrm{-}{\mult}_{{\mu}:x_0} X= {2}/{r}$. \\  
Case$(2)$ $cAx/4$. Let $(X, x_0)$ be $(x^2+y^2 +f(z, u^2)=0 
\subset \mathbb{C}^4/ \mathbb{Z}_4(1,3,2,1), 0)$.
We consider that, for $x^2 = - y^2 - f(z, u^2)$,  
\[ \frac{ O_{X, x_0} }
{ {\mu}_{*} O_{Y} ( -  h {E} ) } \ni    
     \begin{array}{@{\,}ll} 
          {y}^t {z}^l {u}^m \textrm{\ for\ } 
          t + l + m < h \textrm{\ and\ }
          3t + 2 l + m \in 4 \mathbb{Z} 
          \textrm{\ or\ } \\
          x {y}^t {z}^l {u}^m \textrm{\ for\ } 
          1 + t + l + m < h \textrm{\ and\ }
          1 + 3 t + 2 l + m \in 4 \mathbb{Z}.
     \end{array} \]
We fixed $j_1$ and $j_2$ for $0 \leq j_1, j_2 < 4$ 
such that  
$t = 4 n_{j_1} + {j_1}(n_{j_1} \geq 0)$ and  
$l = 4 n_{j_2} + {j_2}(n_{j_2} \geq 0)$.
Then  
we have $m= 4 n - \overline{3 {j_1} + 2 {j_2}}$ or 
$m= 4 n' - \overline{3 {j_1} + 2 {j_2}} - 1$  
Then $0 \leq 
4 (n_{j_1} + n_{j_2} + n ) + (j_1 + j_2 - 
\overline{3 {j_1}+ 2 {j_2}}), 4 (n_{j_1} + n_{j_2} + n' ) + (j_1 + j_2 
- \overline{3 {j_1}+ 2 {j_2}})- 1 \leq h - 1$. 
We have $- 4 <  j_1 + j_2  - \overline{3 {j_1} + 2 {j_2}} < 8$. 
Then we have 
\[ 2 \cdot 4^2 \frac{( {}_{\lfloor} {(h - 9)/ 4}_{\rfloor} + 3 )!}
{( {}_{\lfloor} {(h - 9)/ 4}_{\rfloor})! 3!}
< \dim \frac{O_{X,x_0}}{{\mu}_{*}O_{Y}(-h {E})} 
< 2 \cdot 4^2 \frac{( {}_{\lfloor} {(h + 3)/ 4}_{\rfloor} + 3 )!} 
{({}_{\lfloor} {(h + 3)/ 4}_{\rfloor})! 3!}. \]
Hence, we have ${\wgh}\textrm{-}{\mult}_{\mu:x_0} X= 1/ 2$. \\
Case $(3)$ $cAx/2$. Let $(X, x_0)$ be 
$(x^2 + y^2 + f(z, u) = 0 \subset \mathbb{C}^4/\mathbb{Z}_2(1, 2, 1, 1), 0)$ 
where $\ord(f) \geq 2$ for wt$(z, u)=(1, 1)$.\\
Then we consider that, for $x^2 = - y^2 - f(z, u)$, 
\[ \frac{{\mu}_{*} O_{Y} (-h {E})}
{ {\mu}_{*} O_{Y} (-( h + 1 ){E})}
\ni \begin{array}{@{\,}ll} 
       {y}^t {z}^l {u}^m \textrm{\ for\ } 
       t + l + m = h \textrm{\ and\ } 
       2 t + l + m \in 2 \mathbb{Z} \textrm{\ or\ } \\
       x {y}^t {z}^l {u}^m \textrm{\ for\ } 
       1 + t + l + m = h \textrm{\ and\ }
       1 + 2 t + l + m \in 2 \mathbb{Z} . 
    \end{array} \]
Then we have 
\[ \dim \frac{{\mu}_{*} O_{Y} (-h E)}
{{\mu}_{*} O_{Y} (-( h + 1)E)} 
= \sum_{t=0}^{h} \{ (h - t ) + 1 \}
= h^2 / 2 + 3 h / 2 +1. \]
Hence, we have ${\wgh}\textrm{-}{\mult}_{\mu:x_0} X= 1$.  
The proofs of cases $(4)$ and $(6)$ are the same as 
the proof of case $(3)$. \\
Case$(5)$ $cD/3$. 
Let $(X,x_0)$ be ${u}^2 + {x}^3 + {y}^3 +{z}^3 = 0$ or 
$\{{u}^2 + {x}^3 + y {z}^2 + f(x, y, z) = 0$, 
$f \in (x, y, z)^4$\} or
\{ ${u}^2 + {x}^3 + {y}^3 + f(x, y, z) = 0$, $f \in (x, y, z)^4$\} 
$\subset \mathbb{C}^4/ \mathbb{Z}_3( 1, 2, 2, 0 )$. \\  
We consider that, for $u^2 =  - h (x,y,z)$, 
\[\frac{O_{X, x_0}}{{\mu}_{*}O_{Y}(- h E)} 
\ni \begin{array}{@{\,}ll}
      {x}^t {y}^s {z}^m \textrm{\ for\ } 
      t + s + m < h \textrm{\ and\ }  
      t + 2s + 2 m \in 3 {\mathbb Z} 
      \textrm{\ or\ } \\ 
      {x}^t {y}^s {z}^m u  \textrm{\ for\ } 
      t + s + m + 1 < h \textrm{\ and\ }
      t + 2s + 2 m \in 3 {\mathbb Z}.
    \end{array} \]    
As in case $(2)$, we have ${\wgh}\textrm{-}\mult_{\mu:x_0} X= {2}/{3}$.
\end{proof}
Fourth we calculate the weighted multiplicities of quotient of rational 
double points under the action $\mathbb{Z}_r (1, a, b)$ for 
$(r, a)=1$ and $0 \leq b < r$.  
\begin{Lemm}
Let $(X, x_0)$ be a $3$ folds nonhypersurface and not quotient 
terminal singular point of $\ind_{x_0} X = r >0$ 
over $\mathbb{C}$ and $\mu:Y \rightarrow X$ 
the weighted blow up with the weights 
${\wght}(x, y, z, u) = (1, 1, 1, 1)$  
with the exceptional divisor ${E}$ of ${\mu}$
such that $K_Y = {\mu}^{*} K_X + E$.    
Let 
$C \cong \mathbb{C} / {\mathbb Z}_r (1)$ and 
$S \cong $
$(1)x y + z^{n}=0 \subset \mathbb{C}^3 / {\mathbb Z}_r (1, a, b)$, or 
$(2)x y + z^{n}=0 \subset \mathbb{C}^3 / {\mathbb Z}_r (b, 1, a)$, 
$(3)x^2 + f(y, z) =0 \subset \mathbb{C}^3 / {\mathbb Z}_r (1, a, b)$, 
$(4)x^2 + f(y, z) =0 \subset \mathbb{C}^3 / {\mathbb Z}_r (b, 1, a)$ 
for $(r, a)=1$ and $0 \leq b < r$.   
Then 
\[{\wgh}\textrm{-}\mult_{\mu:x_0} C =1 / r 
\textrm{\ and\ }
{\wgh}\textrm{-}\mult_{\mu:x_0} S =2 / r.\]
\end{Lemm}
\begin{proof}  
Let $\bar{C}$ the strict transform of $C$, 
$\bar{S}$ the strict transform of $S$, and 
$e$ such that $a e \equiv 1 ( \textrm{\ mod\ } r)$. \\
In $C$ case, we consider that, 
\[
\frac{ O_{C, x_0} }{ {\mu}_{*}O_{\bar{C}}(-h E |_{\bar{C}})} 
\ni x^s\textrm{\ for\ } 
h  > s \textrm{\ and\ } s \in r \mathbb{Z}. \]    
As in Theorem 2.5, 
${\wgh}\textrm{-}{\mult}_{{\mu}:x_0} C= 1 / r$. \\
In the $(1)$ case, we consider that, 
\[
\frac{ O_{S, x_0} }{ {\mu}_{*}O_{\bar{S}}(-h E |_{\bar{S}})} 
\ni \begin{array}{@{\,}ll} 
       x^s z^u \textrm{\ for\ } 
       h  > s + u \textrm{\ and\ } 
       s + b u \in r \mathbb{Z}
        \textrm{\ or\ } \\
       y^t z^u \textrm{\ for\ } 
       h > t + u \textrm{\ and\ } 
       t + \overline{e b} u \in r \mathbb{Z}. 
    \end{array} \]    
As in Theorem 2.5, ${\wgh}\textrm{-}\mult_{\mu:x_0} S = 2/ r$. \\ 
In the $(2)$ case, we consider that 
\[
\frac{ O_{S, x_0} }{ {\mu}_{*}O_{\bar{S}}(-h E |_{\bar{S}})} 
\ni \begin{array}{@{\,}ll} 
       x^s z^u \textrm{\ for\ } 
       h  > s + u \textrm{\ and\ }
       \overline{e b} s + u \in r \mathbb{Z}
        \textrm{\ or\ } \\
       y^t z^u \textrm{\ for\ } 
       h > t + u \textrm{\ and\ }
       {e}t + u \in r \mathbb{Z}. 
    \end{array} \]    
As in Theorem 2.5, 
${\wgh}\textrm{-}\mult_{\mu:x_0} S = 2 / r$. \\
In the $(3)$ case, we consider that 
\[ \frac{ O_{S, x_0} }{{\mu}_{*}O_{\bar{S}}(-h E |_{\bar{S}})} 
\ni \begin{array}{@{\,}ll} 
       y^t z^u \textrm{\ for\ } 
       h >t + u \textrm{\ and\ }
       t + \overline{e b} u \in r \mathbb{Z}
        \textrm{\ or\ } \\
        x y^{t} z^{u} \textrm{\ for\ } 
       h >1 + t + u \textrm{\ and\ }
        e + t +\overline{e b} u \in r \mathbb{Z}.
    \end{array} \]   
As $S_3$ case in Lemma 2.6, 
${\wgh}\textrm{-}\mult_{\mu:x_0} S= 2/r$. \\
In the $(4)$ case, we consider that, 
\[
\frac{ O_{S, x_0} }{{\mu}_{*}O_{\bar{S}}(-h E |_{\bar{S}})} 
\ni \begin{array}{@{\,}ll} 
       y^t z^u \textrm{\ for\ } 
       h >t + u \textrm{\ and\ }
       t + a u \in r \mathbb{Z}
        \textrm{\ or\ } \\
        x y^{t} z^{u} \textrm{\ for\ } 
       h >1 + t + u \textrm{\ and\ }
        b + t +a u \in r \mathbb{Z}
    \end{array} \]   
As $S_3$ case in Lemma 2.6, ${\wgh}\textrm{-}\mult_{\mu:y_0} S= 2 /r$. 
\end{proof}
\section{General methods for freeness of adjoint linear systems}
\hspace{1.5 em}
We can construct divisors which have high weighted order 
at a given point from the following: 
\begin{Lemm} 
Let $X$ be a normal and complete variety of dimension $n$,
$L$ a nef and big ${\mathbb Q}$-Cartier divisor, 
$x_0 \in X$ a point, and  $t$,$t_0$ a rational number
such that $t>t_0>0$. We assume that $\mu: Y \rightarrow  X$ is 
the weighted blow up of $X$ at $x_0$ 
with the exceptional divisor $E$ of $\mu$.
Let $W \subset X$ be a subvariety of dimension $p$ such that 
$W$ is normal at $x_0$. 
Then there exists an effective  ${\mathbb Q}$-Cartier divisor $D_W$ such that
$ D_W \sim_{\mathbb Q} t L|_{W} $ and
\[ {\wgh}\textrm{-}{\ord}_{\mu:x_0} D_W \geq (t_0+\epsilon )
 \sqrt[p]{ \frac{L^p W}{ {\wgh}\textrm{-}{\mult}_{\mu:x_0 }W}} \]
which is a rational number for $0 \leq \epsilon 
\ll \sqrt[p]{ {\wgh}\textrm{-}{\mult}_{\mu:x_0} W/ L^p W}$.
\end{Lemm}
\begin{proof}
We change the multiplicity of subvariety at the point 
with the weighted multiplicity of subvariety at the point for $\mu$ 
in [K1 2.1 ( cf [Ka5 2.1])].
The proof is the same as [K1 2.1 ( cf [Ka5 2.1])].
\end{proof}
The following proposition is the key of 
the proofs of our results of freeness: 
\begin{prop}[{[K1, 2.2] cf [Ka5, 2.3]}]  
Let $X$ be a normal projective 
variety of dimension $n$,
$x_0 \in X$ a $KLT$ point,  
and $L$ an ample ${\mathbb Q}$-Cartier divisor such that 
$K_X + L$ is Cartier at $x_0$.
Assume that there exists an effective ${\mathbb Q}$-Cartier divisor $D$
which satisfies the following conditions: \\
$(1)$ $D \sim_{\mathbb Q} t L $ for a rational number $ t<1$,\\
$(2)$ $(X,D)$ is $LC$ at $x_0$,\\
$(3)$ $\{x_0\} \in CLC(X,D)$.\\
Then  $| K_X + L | $ is free at $x_0$.
\end{prop}
We generalize [Ka5, 2.2]
in which ${\mathbb Q}$-divisors have integral order at $x_0$.
For our purpose, we need to treat ${\mathbb Q}$-Cartier divisors 
of fractional weighted orders at $x_0$.
\begin{Lemm}
Let $X$ be a normal projective variety of dimension $n$,
$x_0 \in \mathbb{C}^n/{\mathbb{Z}_r(1,a_1, \cdots a_{n-1})}$
such that an integer $r > 0$, $(r, a_1)=1$, and 
$0 \leq a_j < r$ for $1 \leq j \leq n-1$, 
$L$ an ample $\mathbb{Q}$-Cartier divisor such that 
$K_X + L$ is Cartier at $x_0$. 
Let $l:= \min \{ i\ |\ a_{j}i \equiv i ( \textrm{\ mod\ r } )
(1 \leq j \leq n-1 ) 
\textrm{\ for\ } 0< i \leq r \}$. 
Let ${\mu}: Y \rightarrow X$ be 
the weighted blow up of $X$ at $x_0$ such that 
$wt(x_0, x_1, \cdots, x_{n-1})= (l/r, l/r, \cdots, l/r)$
with the exceptional divisor ${E}$ of ${\mu}$
such that $K_{Y} ={\mu}^{*} K_X + (n l / r -1)E$. 
Let $W$ be a prime divisor with 
${\wgh}\textrm{-}{\ord}_{\mu:x_0} W = {l d}/{r} \geq {l}/{r}$ 
for an integer $d$,
and $e$, $k$ positive rational numbers such that $de \leq 1$ and
$({k l}/{r})^n < {L^n} {l}^n / r^{n-1}$.
Assume that there exists an effective $\mathbb{Q}$-divisor $D$
such that  $D \sim_{\mathbb Q} L$ and 
${\wgh}\textrm{-}{\ord}_{\mu:x_0} D \geq {l k}/{r}$, 
and moreover that $D \geq e k W$ for any such $D$.
Then there exists a real number $\lambda$ with $0 \leq \lambda < 1$
and $\lambda \leq max \{ 1- d e ,(den)^{-{1}/{(n-1)}} \}$
which satisfies the following condition:
if $k'$ is a positive rational number such that $k' > k$ and
\[ (\lambda \frac{k}{r} )^n + ( \frac{1-de-\lambda }{ 1-\lambda } )^{n-1}
\{ (\frac{k'}{r} + \frac{\lambda de }{1-de-\lambda}\frac{k}{r} )^n
-(\frac{ \lambda k}{r} + 
\frac{\lambda de}{1-de-\lambda}\frac{k}{r} )^n \} 
< \frac{L^n}{r^{n-1}},\]
then there exists an effective ${\mathbb Q}$-divisor $D$ such that 
$D \sim_{\mathbb Q} L$ and ${\wgh}\textrm{-}{\ord}_{\mu:x_0} D 
\geq {l k'}/{r}$.
(If $\lambda=1-de$,  
then the left hand side
of the above inequality should be taken as a limit.)
\end{Lemm}
\begin{proof}
(cf [Ka5, 2.2])
We have ${\wgh}\textrm{-}\mult_{{\mu}:x_0} X 
= r^{n-1} / {l}^n$ by Theorem 2.5.
Let $\bar{k}$ = sup \{$q$; there exists an effective ${\mathbb
Q}$-divisor $D$ such that $D \sim_{\mathbb Q} L$ 
and ${\wgh}\textrm{-}{\ord}_{{\mu}: x_0} D={q} l / {r}$ \}.
Let us define a function  $\phi(q)$ for $q \in {\mathbb Q}$
with  $0 \leq q < \bar{k}$ to be the largest real number such that
$D \geq \phi(q) W$ whenever $D \geq 0$, $D \sim_{\mathbb Q} L$ and 
${\wgh}\textrm{-}{\ord}_{{\mu}: x_0} D = {q l}/{r}$.
Then $\phi$ is a convex function [Ka5, 2.2 l25 -- l28 ]. 
Since $\phi(k) \geq ek$, there exists a real number $\lambda$ such
that $0 \leq \lambda < 1$ 
and $\phi(q) \geq {e(q- \lambda k)}/({1-\lambda})$ for any $q$.

Let $m$ be a large and sufficiently divisible integer and 
$\nu : H^0 (X,mL)\rightarrow O_{X,{x_0}}(mL)
\cong O_{X,{x_0}}$ the evaluation homomorphism.
We consider subspaces 
$V_i ={\nu}^{-1}({{\mu}_{*}O_{Y}(-i E)})$
of $H^0 (X,mL)$ for integers $i$ 
such that ${\lambda k m l}/{r} \leq i \leq {k'm l}/{r}$.
First, we have 
\[ \mathrm{dim} V_{\lceil {\lambda km l }/{r} \rceil}
\geq  \mathrm{dim} H^0 (X,mL) - 
\frac{r^{n-1}}{l^n} \frac{ ( \lambda k m l/ r)^n}{n!} 
+ \mathrm{lower\ terms\ in\ } m. \]  
Let $D \in | mL |$ be a member corresponding to 
$h \in V_i$ for some $i$.
Since we have $D \geq \phi(ir/ l m) m W$,
the number of conditions in order for $h \in V_{i+1}$ is at most 
( the number of homogeneous polynomials of order 
$i-\phi(ir/ l m)m ({d l}/{r})$ in $n$ variables )
$\times$ ${\wgh}\textrm{-}\mult_{{\mu}:x_0} X$, i.e.,
\[ \frac{r^{n-1}}{{l}^n} 
\frac{ ( i- \phi(ir/ l m ) m d l/ r )^{n-1}}{{n-1}!}+
\mathrm{lower\ terms\ in\ } m. \]
Therefore, we have $k' < \bar{k}$ because
\[ \frac{r^{n-1}}{{l}^n} \{ \frac{ (\lambda k m l /r )^n}{n!} +
\sum_{i=\lceil {\lambda km l/r}\rceil}^
{k'm l/r - 1}
\frac{ ( i - \frac{ e ( r i /l m -\lambda k)}
{(1-\lambda)} m {d}l/{r})^{n-1}}{{n-1}!} \} +
\mathrm{lower\ terms\ in\ } m \]
\[= \frac{r^{n-1}}{l^n} \{ \frac{ (\lambda k m l/ r )^n }{n!} +
\frac{ ( \frac{1-de-\lambda}{1-\lambda} )^2 
\{ ( k'l/ r + \frac{ \lambda d e}{1-de-\lambda}{k l}/{r})^n -
( \lambda k l / r +
\frac{ \lambda d e}{1-de-\lambda } {k l}/{r})^n \} m^n }{n!} \}\]
\[+ \mathrm{lower\ terms\ in\ } m 
< \frac{m^n L^n}{n!} + \mathrm{lower\ terms\ in\ } m.\]
By [Ka5 2.2 last part], we have that $\lambda \leq {\max}
\{ 1-de,(nde)^{-{1}/(n-1)} \}$.
\end{proof}
\section{Main Theorem}
\hspace{1.5 em}  
First we consider some isolated quotient singular points. 
\begin{thm}
Let $X$ be a normal projective variety of dimension $3$,
$x_0 \in X$ a quotient singular point 
of type $({1}/{r},{a}/{r},{b}/{r})$
such that an integer $r > 0$, $(a, r)=1$ and $(b,r)=1$, 
and $L$ an ample ${\mathbb Q}$-Cartier divisor such that 
$K_X + L$ is Cartier at $x_0$.
Assume that there are positive numbers ${\sigma}_p$ for $p = 1,2,3$
which satisfy the following conditions: \\
$(1)$ $\sqrt[p]{ L^p \cdot W} \geq {\sigma}_p$ 
for any subvariety $W$ of dimension $p$ which contains $x_0$,\\
$(2)$ ${\sigma}_1 \geq {3}/{r}$,
${\sigma}_2 \geq {3}/{\sqrt{r}}$, and 
${\sigma}_3 >  {3}/{\sqrt[3]{r}}$.\\
Then $| K_X + L |$ is free at $x_0$.
\end{thm}
\begin{proof}
Let $l := \min \{ i\ |\ ai \equiv bi \equiv i ( \textrm{\ mod\ } r)  
\textrm{\ for\ } 0< i \leq r \}$. 
Let $\mu: Y \rightarrow X$ be 
the weighted blow up of $X$ at $x_0$ such that 
$wt(x, y, z)= (l/r, l/r, l/r)$
with the exceptional divisor ${E}$ of $\mu$
such that $K_{Y} ={\mu}^{*} K_X + (3 l / r -1)E$. 
Hence ${\wgh}\textrm{-}{\mult}_{\mu:x_0} X = r^2/ {l}^3$ by Theorem 2.5. \\
Step 0. Let $t$ be a rational number such that  
$t >  ({3}/{\sqrt[3]{r}}) / \sqrt[3]{({L}^3)} $.
Since ${\sigma}_3 > {3}/{\sqrt[3]{r}}$, we can take $t < 1$.
Let $t_0$ be a rational number such that 
$ t_0 = ({3}/{\sqrt[3]{r}}) / \sqrt[3]{{L}^3}-\epsilon$
for $0 \leq \epsilon \ll \sqrt[3]{r^2/{l}^3 L^3}$.
By Lemma 3.1, there exists an effective ${\mathbb Q}$-Cartier divisor
$D$ such that  $D \sim_{\mathbb Q} t L$ and 
${\wgh}\textrm{-}{\ord}_{{\mu}:x_0} D 
\geq (t_0+\epsilon)\sqrt[3]{(L^3){l}^3/r^2}$.
Hence ${\wgh}\textrm{-}{\ord}_{\mu:x_0}D = 3 {l}/r$.
 
Let $c$ be the log canonical threshold of $ ( X,D )$ at $x_0$:
\[ c =\sup{ \{t \in {\mathbb Q}; \mbox{ $K_X +t D$ is $LC$ at $x_0$ } \} }.\]
Then $c \leq 1$.
Let $W$ be the minimal element of $CLC(X,x_0,cD)$.
If $W=\{x_0\}$, then $| K_X+L |$ is free at $x_0$ by Proposition 3.2, 
since $c t <1$.\\
Step 1. We consider the case in which $W =C$ is a curve.
By Proposition 1.3, $C$ is normal at $x_0$, i.e., smooth at $x_0$. 
Let $U$ be a neighborhood at $x_0$.
Then $C |_{U} \cong \mathbb{C}/\mathbb{Z}_r(1)$.
Hence we have ${\wgh}\textrm{-}\mult_{\mu:x_0} C = 1/ l$ 
by Lemma 2.6. 
Since $\sigma_1 \geq 3/r$,
there exists a rational number $t'$
with $c t +(1-c)<t'<1 $ and an effective 
${\mathbb Q}$-Cartier divisor $D_C' $
on $C$ such that  $D_C' \sim_{\mathbb Q}(t'-c t )L|_C$ and 
${\wgh}\textrm{-}{\ord}_{\mu:x_0}D_C' = 3 l (1-c) / r$.
As in [Ka5, 3.1 Step1], there exists an effective ${\mathbb Q}$-Cartier 
divisor $D'$ on $X$ such that $D'\sim_{\mathbb Q}(t'-c t )L$ 
and $ D'|_C=D_C'$.
Let $D_1'$ be a general effective ${\mathbb Q}$-Cartier divisor 
on an affine neighborhood $U_1$ of $ x_0$ in $ X$ such that 
$D_1'|_{C \cap U_1} = D_C'|_{C \cap U_1}$ and 
${\wgh}\textrm{-}{\ord}_{\mu:x_0}D_1'=3l(1-c) / r$.
Then we have  ${\wgh}\textrm{-}{\ord}_{\mu:x_0}(cD+D_1')=3 l/r$, 
hence $\{x_0\} \in CLC(U_1,cD+D_1')$.
Let 
\[ c' =\sup{\{ t \in {\mathbb Q};
\mbox{ $K_X +(cD+tD_1')$ is $LC$ at $x_0$ } \} }. \]
Since $D_1'$ is chosen to be general, we have $c' > 0$.
We have an element $W'$ such that $W'\in CLC(X,x_0,cD+c'D_1')$ and
$W' \not\supseteq C$. By Proposition 1.3, $CLC(X,x_0,cD+c'D_1')$ has
an element which is properly contained in $C$.
By Proposition 1.4, 
we conclude that $(X, c D + c' D')$ is $LC$ at $ x_0$, and
$CLC(X, x_0, c D + c' D')$ has
an element which is properly contained in $ C$, i.e., $\{x_0\}$.\\
Step 2. We consider the case in which $W=S$ is a surface. 
Let $U$ be a neighborhood at $x_0$ 
and $h: \widetilde{U} \longrightarrow U$ 
its closure of local universal cover, \\
i.e., $U={\mathbb C}^3/{\mathbb Z}_r(1, a, b)$
and $\widetilde{U} ={\mathbb C}^3$.
Let $\widetilde{S}$ be $h^{*} (S |_{U} )$ and $x_1$ $h^{-1} (x_0)$. 
Since $x_0 \in S$ is a $KLT$ point [Ka5, 1.7] and 
$h$ is ramified only over $x_0$, $\widetilde{S} \subset {\mathbb C}^3$ is 
also irreducible and $x_1 \in \widetilde{S}$ is at most $KLT$.
Moreover $\widetilde{S} \subset {\mathbb C}^3$ is a hypersurface, 
$\widetilde{S}$ is also Gorenstein. Therefore, 
$x_1 \in \widetilde{S}$ is a rational Gorenstein point, that is, 
$x_1 \in \widetilde{S}$ is a smooth point or a rational double point.
If $x_1$ is a smooth point of $\widetilde{S}$, 
then $S |_{U} \cong {\mathbb C}^2/{\mathbb Z}_r(1, a')$ for $(a', r)=1$.
By Lemma 2.6, 
we have ${\wgh}\textrm{-}\mult_{\mu:x_0} S= r / l^2$. 
If $x_1$ is a rational double point of $\widetilde{S}$, 
then $\widetilde{S}$ is A type or D type or E type, i.e, 
$\widetilde{S} \cong x^2 + f(y,z) =0 \textrm{\ or\ } xy + z^{n+1}=0( n >0)$ 
and $\widetilde{S}$ is invariant by the action of $\mathbb{Z}_r(1,a',b')$ 
for $(r, a') =1$ and $(r, b')=1$. 
We have ${\wgh}\textrm{-}\mult_{\mu:x_0} S= 2 r / l^2$ by Lemma 2.6. 
Let $d:={\wgh}\textrm{-}{\mult}_{\mu:x_0} S ( l^2 / r ) = 1$ or $2$. \\
Step 2-1.
We assume first that $d = 1$.
As in Step 1, 
there exists a rational number $t'$
with $c t + (1 - c) < t' < 1$.
By Lemma 3.1 and ${\sigma}_2 \geq 3/ \sqrt{r}$,  
there exist an effective ${\mathbb Q}$-Cartier divisor $D'$ on $X$ 
and a positive number $c'$ such that 
$D' \sim_{\mathbb Q} (t'-c t ) L$, 
${\wgh}\textrm{-}{\ord}_{\mu:x_0}D'|_{S} =(3 l/r)(1 - c)$, 
$(X, c D + c' D')$ is $LC$ at $ x_0$, and
that the minimal element $W'$ of $CLC(X,x_0,cD+c'D')$ is 
properly contained in $S$. Thus we have the theorem when $W' = \{x_0\}$.  
We consider the case in which $W'=C$ is a curve.
Since $t, t' < 1 $, we have $c t + c'(t'-c t)+(1-c)(1-c') <1$.
As in Step 1, we take a rational number $t''$, 
an effective ${\mathbb Q}$-Cartier divisor $D''$ on $X$, 
and a positive number $c''$ such that  
$c t + c'(t'-c t)+(1-c)(1-c') <t''<1$,  
$D''\sim_{\mathbb Q}(t''-c t-c'(t'-c t))L$, 
${\wgh}\textrm{-}{\ord}_{x_0}D''|_C=(3 l/ r)(1-c)(1-c')$,
$(X, c D + c' D' + c'' D'')$ is $LC$ at $x_0$ and 
that the minimal element $W''$ of $CLC(X,x_0,cD+c'D'+c''D'')$ is
properly contained in $C$, i.e., $ \{ x_0 \} \in
CLC(X,x_0,cD+c'D'+c''D'')$. \\
Step 2-2.
We assume that $d = 2$.
As in Step 2-1,
we take a rational number $t'$ with $c t +\sqrt{2}(1-c) < t'$
and  an effective ${\mathbb Q}$-Cartier divisor $D'$ on $X$ with
$D' \sim_{\mathbb Q}(t'-c t)L$ and 
${\wgh}\textrm{-}{\ord}_{\mu:x_0}D'|_S =(3 l / r)(1-c)$.
Here we need the factor $\sqrt{2}$
because $S$ has the weighted multiplicity $2r / l^2$ at $x_0$.
Then we take $0 < c' \leq 1$ such that 
$(X, c D +c' D')$ is $LC$ and $CLC(X, x_0, c D+c' D')$ has 
an element which is properly contained in $S$. 

We shall prove that we may assume $c t +\sqrt{2}(1 - c) < 1$.
Then we can take $t' < 1$ as in Step 2-1,
and the rest of the proof is the same.
For this purpose, we apply Lemma 3.3.
In argument of Steps 0 through 2-1, the number $t$ was chosen
under the only condition that $t<1$.
So we can take $t=1-\epsilon_1$, where the $\epsilon_n$ for 
$n = 1,2,...$ will stand for very small positive rational numbers.
Then $kl/r= 3l/r(1-\epsilon_1) 
=3l/r +\epsilon_2$ and $e=1/3c$.
This means the following:
for any effective $D \sim_{\mathbb Q} tL$, 
if ${\wgh}\textrm{-}{\ord}_{\mu:x_0} D \geq 3 l / r$,
then $cD \geq S$.
We look for $k'> 6/(3-\sqrt{2})$
so that there exists an effective ${\mathbb Q}$-Cartier divisor
$D \sim_{\mathbb Q} tL$ with $t < (3-\sqrt{2})/2$
and ${\wgh}\textrm{-}{\ord}_{\mu:x_0} D \geq 3 l/r$.
The equation for $k'$ becomes
\[ {\lambda}^{3} +( \frac{1-2e-\lambda}{1-\lambda})^2
\{(\frac{2}{3-\sqrt{2}}+\frac{2 \lambda e}{1-2e-\lambda})^{3}
-(\lambda + \frac{2 \lambda e}{1-2e-\lambda})^{3} \} < 1.  \]
We have $\lambda \leq 1/ \sqrt{6e}$, 
$1/3 \leq e \leq 1/2$, and 
in particular, $0 \leq \lambda \leq 1/ \sqrt{2}$.
By [Ka5, 3.1 Step2-2],
we obtain a desired $D$, and can choose a new $t$
such that $t < (3-\sqrt{2})/2$.
Then we repeat the preceding argument from Step 0.
If we arrive at Step 2-2 again, then we have 
$2/3 \leq c \leq 1$ and $ct +\sqrt{2}(1-c) < 1$.  
\end{proof}
The following shows that the conditions in Theorem 4.1 is best possible. 
\begin{Exam}
Let $X=\mathbb{P}(1,1,1,r)$ and $x_0=(0:0:0:1)$.
Then $x_0$ is a quotient singular point of type 
$\mathbb{C}^3 / \mathbb{Z}_r(1,1,1)$ and
$K_X=\mathcal{O}(-r-3)$.
If $K_X+L$ is Cartier at $x_0$ and $L$ is effective, we have 
$L=\mathcal{O}(rk+3)(k \in \mathbb{Z}, rk+3 \geq 0)$.
If $L=\mathcal{O}(3)$, $S=\mathbb{P}(1,1,r)$, and 
$C=\mathbb{P}(1,r)$, then $|K_X+L|$ is not free at $x_0$ and 
we have $L^3=27/r$, $L^2 S = 9 / r$, and $L C = 3 / r$.
\end{Exam}
We have the following  
that $K_X +L$ is not free at a quotient terminal singular point 
for $L^3 > 27 /r$ but $L C < 3 /r$.
\begin{Exam}
Let $X=\mathbb{P}(1,a,r-a,r)$ for $r > 2a$ and $x_0=(0:0:0:1)$.
Then $x_0$ is a quotient singular point of type 
$\mathbb{C}^3 / \mathbb{Z}_r(1,a,r-a)$ and $K_X=\mathcal{O}(-2r-1)$.
If $K_X+L$ is Cartier at $x_0$ and $L$ is effective, we have 
$L=\mathcal{O}(rk+1)(k \in \mathbb{Z}, rk+1 \geq 0)$.
If $L=\mathcal{O}(r+1)$  
and $C= \mathbb{P}(r-a, r)$, 
then $|K_X+L|$ is not free at $x_0$ and 
$L^3=(r+1)^3/ra(r-a) > 27/r$ but $LC = (r+1)/r(r-a) < 3/r$.
\end{Exam} 
We obtain estimates for nonhypersurface and not quotient 
${\mathbb Q}$-factorial terminal singularities. 
\begin{thm}
Let $X$ be a normal projective variety of dimension $3$, 
$x_0 \in X$ a nonhypersurface and not quotient 
$\mathbb{Q}$-factorial terminal singular point of $\ind_{x_0} X =r >1$, and 
$L$ an ample ${\mathbb Q}$-Cartier divisor 
such that $K_X + L$ is Cartier at $x_0$. 
Assume that there are positive numbers ${\sigma}_p$ for $p = 1,2,3$ 
which satisfy the following conditions: \\
$(1)$ $\sqrt[p]{ L^p \cdot W } \geq {\sigma}_p $
for any subvariety $W$ of dimension $p$ which contains $x_0$,\\ 
$(2)$ ${\sigma}_1 \geq 2/r$, ${\sigma}_2 \geq 2\sqrt{2/r}$, 
and ${\sigma}_3 > 2\sqrt[3]{2/r}$.\\
Then $| K_X + L |$ is free at $x_0$.
\end{thm}
\begin{proof}
Let $\mu:Y \rightarrow X$ 
the weighted blow up with the weights 
${\wght}(x, y, z, u) = (1, 1, 1, 1)$  
with the exceptional divisor ${E}$ of ${\mu}$
such that $K_Y = {\mu}^{*} K_X + E$.  
Then 
${\wgh}\textrm{-}{\mult}_{\mu:x_0} X= 2 / r$ by Theorem 2.7. \\
Step 0.    
Since $\sigma_3 > 2\sqrt[3]{2/r}$ and Lemma 3.1, 
as in Step 0 of the proof of Theorem 4.1, 
there exist a rational number $t$ and 
an effective $\mathbb{Q}$-Cartier divisor $ D$
such that 
$0 < t <1$, 
$D \sim_{\mathbb Q} t L $, and 
${\wgh}\textrm{-}{\ord}_{\mu:x_0} D =2$.

Let $c$ be the 
\textit{log canonical threshold} of $ ( X,D )$ at $x_0$:
\[ c=\sup{ \{t \in {\mathbb Q}; \mbox{ $K_X +tD$ is $LC$ at $x_0$ } \} }. \]
Then $c \leq 1$. 
Let $W$ be the minimal element of $CLC(X,x_0,cD)$. 
If $W=\{x_0\}$, then $| K_X+L |$ is free at $x_0$ by Proposition 3.2, 
since $c t <1$.\\ 
Let $U$ be a neighborhood at $x_0$, and 
$h: \widetilde{U} \rightarrow U$ its closure of local universal cover, 
i.e., \\
$U_1 = \widetilde{U_1}/\mathbb{Z}_r(a, r-a, r, 1)$ and 
$\widetilde{U_1} = 
\{x y + f(z, {u}^r)=0 \} \subset \mathbb{C}^4$, \\
$U_2 = \widetilde{U_2}/\mathbb{Z}_4(1, 3, 2, 1)$ and 
$\widetilde{U_2} = 
\{ x^2+y^2+f(z, u^2)=0 \} \subset \mathbb{C}^4$, \\ 
$U_3 = \widetilde{U_3} / \mathbb{Z}_2(1, 2, 1, 1)$ and 
$\widetilde{U_3} = \{ x^2 + y^2 +f(z, u)=0,  
f \in (z, u)^4\mathbb{C}\{ z, u \} \} \subset \mathbb{C}^4 $, \\
$U_4 = \widetilde{U_4} / \mathbb{Z}_2(1, 1, 2, 1)$ and 
$\widetilde{U_4} = \{u^2+z^3+x y z + f(x, y)=0, f \in (x, y)^4\}, 
\textrm{or } \{u^2+x y z + {z}^n + f(x, y)=0, f \in (x, y)^4, n \geq 4\}, 
\{u^2+y^2 z + {z}^n + f(x,y)=0, f\in(x, y)^4, n \geq 3 \}
\subset \mathbb{C}^4$, \\
$U_5= \widetilde{U_5}/ \mathbb{Z}_3(1, 2, 2, 3)$ and 
$\widetilde{U_5}= u^2+x^3+y^3+z^3=0, \textrm{or } 
\{u^2+x^3+y z^2+f(x, y, z)=0, f\in (x, y, z)^4\}, 
\{u^2+x^3+y^3+f(x, y, z)=0, f\in (x, y, z)^4 \}
\subset \mathbb{C}^4$, \\
$U_6= \widetilde{U_6} /\mathbb{Z}_2(2, 1, 1, 1)$ and 
$\widetilde{U_6}= \{ u^2+x^3 +g(y, z) x+ h(y, z)=0, g, h 
\in \mathbb{C}\{y, z\}, g, h \in (y, z)^4 \}
\subset \mathbb{C}^4$. \\
Step 1. We consider the case in which $W =C$ is a curve.
By Proposition 1.3, $C$ is smooth at $x_0$.  
Since 
$\{ x = y = u = 0 \subset U_1 \}$,  
$\{ x = y = u = 0 \subset U_2 \}$, 
$\{ x = z = u = 0 \subset U_3 \}$, 
$\{ x = y = u = 0 \subset U_4 \}$, 
$\{ x = y = z = 0 \subset U_5 \}$, and 
$\{ x = y = z = 0 \subset U_6 \}$ are points, 
then $C |_{U} \cong {\mathbb C}/{\mathbb Z}_r(1)$.
Hence by Lemma 2.8, 
we have ${\wgh}\textrm{-}\mult_{\mu:x_0} C = 1/ r$.  
Since $\sigma_1 \geq 2/r$, 
as in Step 1 of the proof of Theorem 4.1, 
we take a rational number $t'$ 
such that $c t + (1-c)< t' <1$, and 
an effective ${\mathbb Q}$-Cartier divisor $D' $ on $X$ 
and a positive number $c'$ 
such that  $D' \sim_{\mathbb Q}(t'-c t )L$, 
${\wgh}\textrm{-}{\ord}_{x_0}D'|_C =2(1-c)$, 
$(X,cD+c'D')$ is $LC$ at $ x_0$, and
the minimal element $W'$ of $CLC(X,x_0,cD+c'D')$ 
is properly contained in $ C$, i.e., $\{x_0\}$.\\
Step 2. We consider the case in which $W = S$ is a surface.  
Let $\widetilde{S}$ be $h^{*}(S|_U)$, and $x_1$ be $h^{-1} (x_0)$. 
Since $\widetilde{U}$ is ${\mathbb Q}$-factorial 
Gorenstein terminal at $x_1$,
$\widetilde{U}$ is factorial at $x_1$([Ka2, Lemma 5.1]).
In particular, $\widetilde{S}$ is a Cartier divisor at $x_1$.
Therefore, $x_1 \in \widetilde{S}$ is a rational double point.
Then $\widetilde{S}$ is A type or D type or E type, i.e, 
$\widetilde{S} \cong x^2 + f(y,z) =0 \textrm{\ or\ } xy + z^{n+1}=0( n >0)$ 
and $\widetilde{S}$ is invariant by the action of $\mathbb{Z}_r(a',1,b')$ 
for $(r, a') =1$ or $(r, b')=1$. 
By Lemma 2.8, 
We have ${\wgh}\textrm{-}\mult_{\mu:x_0} S= 2 / r$. \\
Step 2-1. As in Step 1 of the proof of Theorem 4.1 , 
there exists a rational number $t'$
with $c t + (1 - c) < t' < 1$.
By Lemma 3.1 and ${\sigma}_2 \geq 2/ \sqrt{r}$,  
there exist an effective ${\mathbb Q}$-Cartier divisor $D'$ on $X$ 
and a positive number $c'$ such that 
$D' \sim_{\mathbb Q} (t'-c t ) L$, 
${\wgh}\textrm{-}{\ord}_{\mu:x_0}D'|_S =2 (1 - c)$, 
$(X, c D + c' D')$ is $LC$ at $ x_0$, and
that the minimal element $W'$ of $CLC(X,x_0,cD+c'D')$ is 
properly contained in $S$. Thus we have the theorem when $W' = \{x_0\}$.  
We consider the case in which $W'=C$ is a curve.
Since $t, t' < 1 $, we have $c t + c'(t'-c t)+(1-c)(1-c') <1$.
As in Step 1, we take a rational number $t''$, 
an effective ${\mathbb Q}$-Cartier divisor $D''$ on $X$
and a positive number $c''$ such that  
$c t + c'(t'-c t)+(1-c)(1-c') <t''<1$,  
$D''\sim_{\mathbb Q}(t''-c t-c'(t'-c t))L$, 
${\wgh}\textrm{-}{\ord}_{x_0}D''|_C=2(1-c)(1-c')$,
$(X, c D + c' D' + c'' D'')$ is $LC$ at $x_0$ and 
that the minimal element $W''$ of $CLC(X,x_0,cD+c'D'+c''D'')$ is
properly contained in $C$, i.e., $ \{ x_0 \} \in
CLC(X,x_0,cD+c'D'+c''D'')$.
\end{proof}
The results for smooth case [Ka5 Theorem 3.1] and
for Gorenstein $\mathbb{Q}$-factorial terminal case [K1 Theorem 3.3]
and Theorem 4.1, 4.4 imply the following corollary 
for $\mathbb{Q}$-factorial terminal singularities on threefolds:
\begin{Coro}
Let $X$ be a normal projective variety of dimension $3$ 
with only $\mathbb{Q}$-factorial terminal singularities, 
$x_0 \in X$ a terminal singular point 
such that $\ind_{x_0} X=r > 0$, 
and $L$ an ample ${\mathbb Q}$-Cartier divisor such that 
$K_X + L$ is Cartier at $x_0$.
Assume that there are positive numbers ${\sigma}_p$ for $p = 1,2,3$
which satisfy the following conditions: \\
$(1)$ $\sqrt[p]{ (L)^p \cdot W} \geq {\sigma}_p$ 
for any subvariety $W$ of dimension $p$ which contains $x_0$,\\
$(2)$ ${\sigma}_1 \geq {3}/{r}$,
${\sigma}_2 \geq {3}/{\sqrt{r}}$, and 
${\sigma}_3 >  {3}/{\sqrt[3]{r}}$.\\
Then $| K_X + L |$ is free at $x_0$.
\end{Coro}
\section{Effective freeness conditions on surfaces 
at log terminal singularities of A type }
\hspace{1.5 em}
We know that, for effective freeness condition,  
the weighted blow up and the weighted multiplicities are more important 
than the fundamental cycle and Kawachi invariant introduced by Sakai [S] 
and more powerful Kawachi invariant. 
We give an alternative proof of a theorem of Langer for
log terminal singular points of only A type [La2]. 
The method of [La2] is to use a reflexive sheaf of rank $2$, 
which is completely different for ours.  
\begin{thm}[{[La2, Theorem 0.2's application]}]
Let $S$ be a normal projective surface,
$x_0 \in S$ log terminal singularities of $A$ type, i.e., 
$\cong \mathbb{C}^2/{\mathbb{Z}_r(a,1)}$ 
such that an integer $r > 0$ and $(r, a) = 1$, 
and $L$ an ample ${\mathbb Q}$-Cartier divisor such that 
$K_S +L$ is Cartier at $x_0$.
Assume that there are positive numbers ${\sigma}_p$ for $p =1,2$
which satisfy the following conditions:\\
$(1)$ $\sqrt[p]{ L^p \cdot W} \geq {\sigma}_p$
for any subvariety $W$ of dimension $p$ which contains $x_0$, \\
$(2)$ ${\sigma}_1 \geq 2 / r$ 
and ${\sigma}_2 > 2 / \sqrt{r}$. \\
Then $| K_S + L |$ is free at $x_0$.
\end{thm}
\begin{proof}
Let $l := \min \{ i\ |\ ai \equiv i ( \textrm{\ mod\ } r ) 
\textrm{\ for\ } 0< i \leq r \}$ and 
${\mu}: \bar{S} \rightarrow S$ be 
the weighted blow up of $S$ at $x_0$ such that $wt(x, y)= (l/r, l/r)$
with the exceptional divisor ${E}$ of ${\mu}$
such that $K_{\bar{S}} ={\mu}^{*} K_S + (2 l / r -1)E$. 
By Theorem 2.5, ${\wgh}\textrm{-}{\mult}_{\mu:x_0} S = r/ l^2$. \\
Step 0. 
Since $\sigma_2 >2 / \sqrt{r}$ and Lemma 3.1, 
as in Step 0 of the proof of Theorem 4.1, 
there exists a rational number $t$ and 
an effective ${\mathbb Q}$-Cartier divisor $D$
such that 
$0 < t < 1$, $D \sim_{\mathbb Q} t L $, and 
${\wgh}\textrm{-}{\ord}_{\mu:x_0} D =2 l / r$. 

Let $c$ be the 
\textit{log canonical threshold} of $ ( S,D )$ at $x_0$:
\[ c=\sup{ \{t \in {\mathbb Q}; \mbox{ $K_S +tD$ is $LC$ at $x_0$ } \} }. \]
Then $c \leq 1$. 
Let $W$ be the minimal element of $CLC(S,x_0,cD)$.
If $W=\{x_0\}$, then $|K_S+L|$ is free at $x_0$ by Proposition 3.2, 
since $c t <1$.\\
Step 1. We consider the case in which $W =C$ is a curve.
By Proposition 1.3, $C$ is smooth at $x_0$. 
Let $U$ be a neighborhood at $x_0$. 
Then $C|_{U} \cong \mathbb{C}/ \mathbb{Z}_r(1)$. 
Hence we have ${\wgh}\textrm{-}\mult_{\mu:x_0} C = 1/ l$. 
Since $\sigma_1 \geq 2 / r$, as in Step 1 of Theorem 4.1, 
we take a rational number $t'$, 
an effective ${\mathbb Q}$-Cartier $D'$ on $X$, 
and a positive number $c'$ 
such that $c t + ( 1 - c ) < t' < 1$, $D' \sim_{\mathbb Q}(t'-ct )L$,  
${\wgh}\textrm{-}{\ord}_{x_0}D'|_C =(2l/r)(1-c)$,
$(S, c D + c' D')$ is $LC$ at $x_0$, and
that the minimal element $W'$ of $CLC(S,x_0,cD+c'D')$ is 
properly contained in $ C$, i.e., $\{x_0\}$.
\end{proof}
The following shows the conditions in Theorem 5.1 is best possible.
\begin{Exam}
Let $S=\mathbb{P}(1, 1, r)$ 
for an integer $r > 0$ and $x_0=(0:0:1)$. Then 
$x_0$ is a quotient singular point of type 
$\mathbb{C}^2 / \mathbb{Z}_r(1, 1)$ 
and $K_S = \mathcal{O}(- r -2)$. 
If $K_S +L$ is Cartier at $x_0$ and $L$ is effective, 
we have $L=\mathcal{O}(rk + 2)
(k \in \mathbb{Z}, rk + 2 \geq 0)$.
If $L=\mathcal{O}(2)$ and $C = \mathbb{P}(1,r)$, 
then $|K_S +L|$ is not free at $x_0$ and 
$L^2 = 4 / r$ and $LC = 2 /r$.
\end{Exam}
We have the following that $K_S +L$ is not free at log 
terminal singularities of A type for $L^2 > 4 / r$ and $LC < 2/r$. 
\begin{Exam}
Let $S=\mathbb{P}(1, a, r)$ such that $1 < a < r$ 
and $x_0=(0:0:1)$. Then 
$x_0$ is a quotient singular point of type 
$\mathbb{C}^2 / \mathbb{Z}_r(1, a)$ 
and $K_S = \mathcal{O}(- r -a -1)$. 
If $K_S +L$ is Cartier at $x_0$ and $L$ is effective, 
we have $L=\mathcal{O}(rk + a + 1)
(k \in \mathbb{Z}, rk + a +1 \geq 0)$.
If $L=\mathcal{O}(a + 1)$ and $C = \mathbb{P}(a,r)$, 
then $|K_S +L|$ is not free at $x_0$ and 
$L^2 =(a+1)^2/ar > 4 / r$ and 
$LC = 2 / a r < 2 /r$.
\end{Exam}

\end{document}